

\baselineskip=14pt
\parskip=10pt
\def\halmos{\hbox{\vrule height0.15cm width0.01cm\vbox{\hrule height
  0.01cm width0.2cm \vskip0.15cm \hrule height 0.01cm width0.2cm}\vrule
  height0.15cm width 0.01cm}}

\magnification=\magstephalf

\def\1{{\overline{1}}}
\def\2{{\overline{2}}}
\parindent=0pt
\overfullrule=0in

\def\frac#1#2{{#1 \over #2}}
\centerline
{\bf The Jackson-Richmond 4CT Constant is EXACTLY 10/27}
\bigskip
\centerline
{\it Shalosh B. EKHAD and Doron ZEILBERGER}

{\bf Abstract}: In their recent claimed computer-free proof of the Four Color Theorem, David Jackson and Bruce Richmond attempted to use sophisticated ``asymptotic analysis''
to explicitly compute a certain number whose positivity (according to them)  implies this famous theorem.
While the jury is still out whether their valiant attempt holds water, we prove, in this modest note, that
this constant equals exactly 10/27. We also point out that their evaluation of this constant
must be erroneous, for two good reasons. Finally, as an encore, we state many similar, but more complicated, results.

David Jackson and Bruce Richmod recently made a brave attempt [JR] to give a {\it human-generated} proof of the Four Color (or as they would spell it, Colour) Theorem.
They relied heavily on  Tutte's seminal paper [T].
An important part was their (correct) theorem stated below (equivalent to Theorem 9 of [JR]).

(We follow their convention of  denoting the coefficient of $x^n$ in the formal power series $f(x)$ by $[x^n]f(x)$.)

{\bf Theorem}: Let $g(x)$ be the formal power series
$$
g(x)= \sum_{n \geq 1} \frac{2\,(4n+1)!}{(n+1)!(3n+2)!}\, x^n \quad,
$$
then there exists a positive constant, let's call it $c$, such that
$$
\lim_{n \rightarrow \infty} \, \frac{[x^n] (g(x))^2}{[x^n] g(x)}  \, = \, c \quad .
$$

However, as we will show below, their determination of that constant $c$ must be erroneous.

It turns out that this constant, $c$, that may be christened the {\bf Jackson-Richmond 4CT constant}, has a nice explicit expression as a good old {\bf rational number}.

{\bf Fact}:  $c=\frac{10}{27}$.

In fact, we can do more than just  asymptotics. We can actually  state the {\bf exact} expression for that ratio, as a {\bf rational function} of $n$,  and sketch {\bf two} proofs.

{\bf Lemma}: For all integers $n \geq 1$, we have:
$$
\frac{[x^n] (g(x))^2}{[x^n] g(x)} \, = \, \frac{10(n-1)(n^2+14n+12)}{3(3n+5)(3n+4)(n+2)} \quad .
$$

{\bf Sketch of First Proof}: the top of the left side is a certain terminating hypergometric sum, and dividing by $[x^n] g(x)$, alias  $\frac{2(4n+1)!}{(n+1)!(3n+2)!}$, is still such a sum,
and the {\bf Zeilberger algorithm} [Z] finds a linear recurrence equation (second-order in fact) satisfied by it. That same recurrence happens to also hold for the
right side, and checking the initial conditions at $n=1$ and $n=2$ concludes the proof. \halmos \quad .

{\bf Sketch of Second Proof}: Below we will exhibit a degree $4$ algebraic equation satisfied by $g(x)$. Hence $g(x)^2$ satisfies another algebraic
equation of degree $4$. From this using the standard tools described below, one can deduce a differential equation, and in turn, a recurrence
satisfied by the coefficients of $g(x)^2$, that happen to match the right side of the lemma times $[x^n]g(x)$. \halmos \quad .

{\bf Critique of  section 2.2 of [JR]}

In Eq. $(13)$ of [JR], they needed to evaluate $A:=g(\frac{27}{256})$, but failed to give its exact value, that happens to be $\frac{5}{27}$. In fact the
formal power series $g(x)$ (that happens to also be a convergent power series for $x \leq \frac{27}{256}$) satisfies the {\bf algebraic equation}:

$$
x \left( {x}^{2}+11\,x-1 \right) + \left( 4\,{x}^{3}+25\,{x}^{2}-14\,x+1 \right) g(x) + x \left( 6\,{x}^{2}+17\,x+3 \right) {g(x)}^{2}+{x}^{2} \left( 4\,x+3 \right) {g(x)}^{3}+{x}^{3}{g(x)}^{4} \, = 0 \quad .
$$

Once guessed, it is easily confirmed using, Salvy and Zimmermann's [SZ] Maple package {\tt gfun}, specifically the
commands {\tt gfun[algeqtodiffeq]} followed by {\tt gfun[diffeqtorec]}, and then verifying that the sequence
$\{2\frac{(4n+1)!}{(n+1)!(3n+2)!}\}$ also satisfies this recurrence.

Plugging-in $x=\frac{27}{256}$, and solving for $g(\frac{27}{256})$, 
gives the exact value $A=\frac{5}{27}$.

In [JR], $c$ is given as $\frac{27}{2} \sqrt{\frac{3}{2}} \cdot A \cdot B$, where 
$B=\frac{16}{27} \sqrt{\frac{3}{2 \pi}}$. Since the latter, thanks to Ferdinand von Lindemann, is {\bf transcendental}, and we just
exhibited both $A$ and $c$ as rational numbers, the [JR] value must be erroneous. Another good reason why it can't be right, as stated (of course they may be some
minor misprints that we were unable to correct), is that in [JR], $c$ evaluates to $1.253754\dots$, and of course while $c$ is positive (as we proved above), it can't be larger than $1$.

{\bf Encore}

The beauty of symbolic computation is that with barely extra effort, we can do much more! Using the Maple package {\tt Jackmond.txt} accompanying this article 
available from 

{\tt https://sites.math.rutgers.edu/\~{}zeilberg/tokhniot/Jackmond.txt}, 

we discovered the next lemma.

{\bf A More general Lemma:} For all integers $n \geq 1$, and $r \geq 2$, define
$$
A_r(n)\,:=\, \frac{[x^n] (g(x))^r}{[x^n] g(x)}  \quad,
$$
and
$$
B_r:= \lim_{n \rightarrow \infty} A_r(n) \quad ,
$$
We have (as above)
$$
A_2(n)=\frac{10(n-1)(n^2+14n+12)}{3(3n+5)(3n+4)(n+2)} \quad .
$$
$$
B_2 = \frac{10}{27} \quad.
$$
We also have
$$
A_3(n)=
\frac{5 \left(n -1\right) \left(n -2\right) \left(5 n^{4}+160 n^{3}+1803 n^{2}+3768 n +2016\right)}{3 \left(3 n +8\right) \left(3 n +5\right) \left(3 n +7\right) \left(3 n +4\right) \left(n +3\right) \left(n +2\right)} \quad;
$$
$$
B_3 =\frac{25}{243} \quad ;
$$
$$
A_4(n)=
\frac{20 \left(n -1\right) \left(n -2\right) \left(n -3\right) \left(25 n^{6}+1350 n^{5}+31495 n^{4}+347406 n^{3}+1211092 n^{2}+1580304 n +665280\right)}{27 \left(3 n +11\right) \left(3 n +8\right) \left(3 n +5\right) \left(3 n +10\right) \left(3 n +7\right) \left(3 n +4\right) \left(n +4\right) \left(n +3\right) \left(n +2\right)}
$$
$$
B_4 =\frac{500}{19683} \quad .
$$
To see $A_r(n)$ for $5 \leq r \leq 11$ look at the output file

{\tt https://sites.math.rutgers.edu/\~{}zeilberg/tokhniot/oJackmond1.txt} \quad .

The list of the $B_r$ for $2 \leq r \leq 11$ is
$$
\left[{\frac{10}{27}}, {\frac{25}{243}}, {\frac{500}{19683}}, {\frac{
3125}{531441}}, {\frac{6250}{4782969}}, {\frac{109375}{387420489}}, {
\frac{625000}{10460353203}}, {\frac{390625}{31381059609}}, {\frac{
19531250}{7625597484987}}, {\frac{107421875}{205891132094649}}\right] \quad,
$$
and in decimals
$$
[ 0.3703703704,  0.1028806584,  0.02540263171,  0.005880238822, 
 0.001306719738,  0.0002823159928,  
$$
$$
0.00005974941647, 
 0.00001244779510,  0.000002561274712,  0.0000005217411450] \quad .
$$

If you want to see data all the way to $r=18$, feel free to look at

{\tt https://sites.math.rutgers.edu/\~{}zeilberg/tokhniot/oJackmond1a.txt} \quad .

Enjoy!

{\bf Moral}: Before posting your next paper publicly, make sure to do some {\it fact-checking} using our beloved silicon servants (soon to become our masters). 
While, {\it who knows?}, perhaps one day, we won't need them
to prove the Four Color Theorem, since either the [JR] attempt, once corrected, would succeed, or some future smart humans would do it, we can still
use computer-kind for the mundane task of checking our calculations, and confirming our logical arguments, {\it before we go public}.

\vfill\eject

{\bf References}

[JR] D.M. Jackson and L.B. Richmond, {\it A non-constructive proof of the Four Colour Theorem},arxiv,
{\tt https://arxiv.org/abs/2212.09835}

[SZ]  Bruno Salvy and Paul Zimmermann, {\it GFUN: a Maple package for the manipulation of generating and holonomic functions in one variable},
ACM Transactions on Mathematical Software {\bf 20} (1994), 163–177. \hfill\break
{\tt https://dl.acm.org/doi/pdf/10.1145/178365.178368} \quad.

[T] W. T. Tutte, {\it A census of planar triangulations}, Canad. J. Math. {\bf 14} (1962), 21-38.

[Z] Doron Zeilberger, {\it The method of creative telescoping},  J. Symbolic Computation {\bf 11} (1991), 195-204. \hfill\break
{\tt https://sites.math.rutgers.edu/\~{}zeilberg/mamarimY/ct1991.pdf} \quad .

\bigskip
\hrule
\bigskip
Shalosh B. Ekhad, c/o D. Zeilberger, Department of Mathematics, Rutgers University (New Brunswick), Hill Center-Busch Campus, 110 Frelinghuysen
Rd., Piscataway, NJ 08854-8019, USA. \hfill\break
Email: {\tt ShaloshBEkhad at gmail dot com}   \quad .
\bigskip
Doron Zeilberger, Department of Mathematics, Rutgers University (New Brunswick), Hill Center-Busch Campus, 110 Frelinghuysen
Rd., Piscataway, NJ 08854-8019, USA. \hfill\break
Email: {\tt DoronZeil at gmail  dot com}   \quad .
\bigskip
{\bf Exclusively published in the Personal Journal of Shalosh B. Ekhad and Doron Zeilberger, and arxiv.org}
\bigskip
{\bf Feb. 8, 2024} \quad .

\end